\documentclass[A4paper,12pt]{book}
\usepackage{amssymb,amsmath}
\usepackage[russian]{babel}
\usepackage{graphicx}
\newcount\exnom\exnom=0

\long\def\exo/{\vspace{0.2cm} \noindent\advance\exnom by1{\bf
{\the\exnom}.}}

\newcommand{\ds}{\displaystyle}
\newcommand{\ep}{\varepsilon}

\def\Im{\mathop{\rm Im}\nolimits}

\begin{document}
\thispagestyle{empty}

\begin{center}\bf \large IDEMPOTENT ELEMENTS OF THE ENDOMORPHISM SEMIRING OF A FINITE CHAIN
\end{center}

\vspace{1mm}

\begin{center} \textbf{Ivan Trendafilov, Dimitrinka Vladeva}
\end{center}

\vspace{5mm}

\begin{quote}\centerline{{ Abstract}}

 {\small Idempotents yield much insight in the structure of finite semigroups and semi\-rings. In this article, we obtain some results on  (multiplicatively) idempotents of the endomorphism semiring of a finite chain. We prove that the set of all idempotents with certain fixed points is a semiring and find its order. We further show that this semiring is an ideal in a well known semiring. The construction of an equivalence relation such that any equivalence class contain just one idempotent is proposed. In our main result we prove that such equivalence class is a semiring and  find his order. We prove that the set of all idempotents with certain jump points is a semiring. }
\end{quote}

\vspace{8mm}

\centerline{{\large 1. \hspace{0.5mm}Introduction}}

\vspace{5mm}

The idempotents play an essential role in the theory of finite semigroups and semirings. It is
well known that in a finite semigroup  some power of each element is an
idempotent, so the idempotents can be taken to be like a generating system of the semigroup or the semiring. For deep results, using idempotents in the representation theory of finite semigroups we refer to [1] and [4].

Let us briefly survey the contents of our paper. After the preliminaries, in the third part we show some facts about fixed points of idempotent endomorphism. The central result here is Theorem 3.9 where we prove that the set of all idempotents with $s$ fixed points $k_1, \ldots, k_s$, ${1 \leq s \leq n-1}$, is a semiring of order  $\ds \prod_{m=1}^{s-1} (k_{m+1} - k_{m})$. Moreover, this semiring is an ideal of the semiring of all endomorphisms having at least $k_1, \ldots, k_s$ as fixed points. In the next part we consider an equivalence relation on some finite semigroup $S$ such that for any $x, y \in S$ follows $x \sim y$ if and only if $x^k = y^m = e$, where $k, m \in \mathbf{N}$ and $e$ is an idempotent of $S$. Then we consider the equivalence classes of the semigroup $S = \left(\widehat{\mathcal{E}}_{\mathcal{C}_n},\cdot\right)$ which is, see [11], one subsemigroup of $\mathcal{PT}_n$. The main result of the paper is Theorem 4.8 where we prove that such equivalence class is a semiring of order
$$C_{k_{1,1}} \prod_{i =1}^{\ell - 1}\left( C_{t_i} C_{s_i}\right)C_{n-1-k_{\ell,\,m_\ell}},$$
where $C_p$ is the $p$ -- th Catalan number.

 In the last part of the paper we prove that the set of idempotent endomorphisms with the same jump points is a semiring.

\newpage

\centerline{{\large 2. \hspace{0.5mm} Preliminaries}}

\vspace{4mm}

We consider some basic definitions and facts concerning finite semigroups and that can be found in any of [1], [2], [4] and [7].
As the terminology  for semirings  is not completely standardized, we say what our conventions are.

\vspace{2mm}

 An algebra $R = (R,+,.)$ with two binary operations $+$ and $\cdot$ on $R$, is called {\emph{semiring}} if:

$\bullet\; (R,+)$ is a commutative semigroup,

$\bullet\; (R,\cdot)$ is a semigroup,

$\bullet\;$ both distributive laws hold
$ x\cdot(y + z) = x\cdot y + x\cdot z$ and $(x + y)\cdot z = x\cdot z + y\cdot z$
for any $x, y, z \in R$.

\vspace{2mm}

 Let $R = (R,+,.)$ be a semiring.
If a neutral element $0$ of  semigroup $(R,+)$ exists and satisfies $0\cdot x = x\cdot 0 = 0$ for all $x \in R$, then it is called \emph{zero}.
If a neutral element $1$ of  semigroup $(R,\cdot)$ exists, it is called \emph{identity element}.

An element $a$ of a semiring $R$ is called \emph{additively (multiplicatively) idempotent} if ${a + a = a}$ $\;(a\cdot a = a)$.
A semiring $R$ is called \emph{additively idempotent} if each of its elements is additively idempotent, see [3].

An element $a$ of a semiring $R$ is called an \emph{additively (multiplicatively) absorbing element} if and only if
$a + x = a\;\;$ $(a\cdot x = x\cdot a = a)\;\;$ for any $x \in R$.
The zero of $R$ is the unique multiplicative absorbing element; of course it does not need to exist.

\vspace{2mm}

 An algebra$\mathcal{M}$ with binary operation $\vee$ such as

$\bullet$ $\; a\vee(b\vee c) = (a\vee b)\vee c$ for any $a, b, c \in \mathcal{M}$;

$\bullet$ $\; a\vee b = b \vee a$ for any $a, b\in \mathcal{M}$;

$\bullet$ $\; a\vee a = a$ for any $a \in \mathcal{M}$.

is called \emph{semilattice (join semilattice)}.

Another term used for $\mathcal{M}$ is a commutative idempotent semigroup -- see [12].
For any $a, b\in \mathcal{M}$ we denote
$ a \leq b \; \iff \; a \vee b = b$.
In this notation, if there is a neutral element in $\mathcal{M}$, it is the least element.

\vspace{2mm}

For a semilattice $\mathcal{M}$ the set $\mathcal{E}_\mathcal{M}$ of the endomorphisms of $\mathcal{M}$ is a semiring
 with respect to the addition and multiplication defined by:

 $\bullet \; h = f + g \; \mbox{when} \; h(x) = f(x)\vee g(x) \; \mbox{for all} \; x \in \mathcal{M}$,

 $\bullet \; h = f\cdot g \; \mbox{when} \; h(x) = f\left(g(x)\right) \; \mbox{for all} \; x \in \mathcal{M}$.

\vspace{2mm}

 This semiring is called the \emph{ endomorphism semiring} of $\mathcal{M}$.
Here all semilattices are finite chains.
Following [9] and [10] we fix a finite chain $\mathcal{C}_n = \; \left(\{0, 1, \ldots, n - 1\}\,,\,\vee\right)\;$ and denote the endomorphism semiring of this chain with $\widehat{\mathcal{E}}_{\mathcal{C}_n}$. We do not assume that $\alpha(0) = 0$ for arbitrary $\alpha \in \widehat{\mathcal{E}}_{\mathcal{C}_n}$. So, there is not a zero in  endomorphism semiring $\widehat{\mathcal{E}}_{\mathcal{C}_n}$. Subsemiring ${\mathcal{E}}_{\mathcal{C}_n} = {\mathcal{E}}_{\mathcal{C}_n}^0$ of $\widehat{\mathcal{E}}_{\mathcal{C}_n}$ consisting of all endomorphisms $\alpha$ with property $\alpha(0) = 0$ has zero and  is considered in [9], [10] and [12].

\vspace{2mm}

If $\alpha \in \widehat{\mathcal{E}}_{\mathcal{C}_n}$ such that $f(k) = i_k$ for any  $k \in \mathcal{C}_n$ we denote $\alpha$ as an ordered $n$--tuple $\wr\,i_0,i_1,i_2, \ldots, i_{n-1}\,\wr$. Note that mappings will be composed
accordingly, although we shall usually give preference to writing mappings on
the right, so that $\alpha \cdot \beta$ means "first $\alpha$, then $\beta$".
The identity $\mathbf{i} = \wr\,0,1, \ldots, n-1\,\wr$ and all constant endomorphisms $\kappa_i = \wr\, i, \ldots, i\,\wr$ are obviously (multiplicatively) idempotents.

The element $a \in \mathcal{C}_n$ satisfying $\alpha(a) = a$ is usually called a fixed point of the endomorphism $\alpha$. Any ${\alpha \in \widehat{\mathcal{E}}_{\mathcal{C}_n}}$ has at most  $n$ fixed points and only the identity $\mathbf{i}$ has just $n$ fixed points.

\vspace{2mm}

For other properties of the endomorphism semiring we refer to [9], [10],  [11] and [12].

In the following sections we use some terms from book [2] having in mind that in [11] we show that some subsemigroups of the partial transformation semigroup are indeed endomorphism semirings.

\vspace{7mm}

\centerline{{\large 3. \hspace{0.5mm} Idempotent endomorphisms and their fixed points}}

\vspace{3mm}

 The set of all idempotents (or idempotent elements) of  semiring $\widehat{\mathcal{E}}_{\mathcal{C}_n}$ is not a semiring, namely endomorphisms $\alpha = \wr\,0, 0, 2\,\wr$ and $\beta = \wr\, 0, 1, 1\,\wr$ are idempotents of  semiring $\widehat{\mathcal{E}}_{\mathcal{C}_3}$ but $\alpha\beta = \wr\,0, 0, 1\,\wr$ is not an idempotent. One trivial, but useful, condition for the idempotents is the following proposition.

\vspace{3mm}

\textbf{Proposition 3.1 } \textsl{The endomorphism $\alpha \in \widehat{\mathcal{E}}_{\mathcal{C}_n}$ is an idempotent if and only if for any $k \in \mathcal{C}_n$, which is not a fixed point of $\alpha$, the image $\alpha(k)$ is a fixed point of $\alpha$.}

\emph{Proof.} Let $\alpha(k) = j$ implies $\alpha(j) = j$. Then $\alpha^2(k) = \alpha(j) = j = \alpha(k)$, what means that $\alpha^2 = \alpha$.

 Conversely, assume that $\alpha^2 = \alpha$ and $\alpha(k) = j$. Now $\alpha^2(k) = \alpha(k)$ for any $k \in {\mathcal{C}_n}$ and  $\alpha^2(k) = \alpha(j)$. Hence $\alpha(j) = \alpha(k)$  and from $\alpha(k) = j$ follows $\alpha(j) = j$.

\vspace{3mm}

Note that in [5] maps $\alpha_k$ are considered   such that  ${\alpha_k(x) = k}$ for all ${x \in \mathcal{M}}$, where $\mathcal{M}$ is not necessarily finite semilattice. These maps are called constant endomorphisms.

For any $\alpha \in \widehat{\mathcal{E}}_{\mathcal{C}_n}$ with one fixed point the cardinality of set $\Im(\alpha)$ is an arbitrary number from 2 to $n-1$.
Indeed, for $\alpha = \wr\, n-2, n-1, \ldots, n-1\,\wr$ with unique fixed point $n-1$ we have $|\Im(\alpha)| = 2$ and for $\beta = \wr\, 1, 2, \ldots, n-2, n -1, n-1\,\wr$ with the same unique fixed point, $|\Im(\beta)| = n-1$. So, the following consequence of Proposition 3.1. is important.

\vspace{3mm}

\textbf{Corollary 3.2 } \textsl{ Endomorphism with only one fixed point is an idempotent if and only if it is a constant.}

\emph{Proof.} Let $k$ is the unique fixed point of endomorphism $\alpha \in \widehat{\mathcal{E}}_{\mathcal{C}_n}$ and $\alpha$ is an idempotent. Let $\ell \neq k$ è $\alpha(\ell) = m$. Using Proposition 3.1 we have that $\alpha(m) = m$ and from the uniqueness of the fixed point follows $m = k$. So, for any $\ell \neq k$ we prove that $\alpha(\ell) = k$, i.e. $\alpha$ is a constant endomorphism.
Conversely, every constant endomorphism is an idempotent.

\vspace{3mm}

By similar reasonings we prove:

\vspace{3mm}

\textbf{Corollary 3.3 } \textsl{Endomorphism $\alpha \in \widehat{\mathcal{E}}_{\mathcal{C}_n}$ with $s$ fixed points $k_1, \ldots, k_s$, $1 \leq s \leq n-1$, is an idempotent if and only if  $\,\Im(\alpha) = \{k_1, \ldots, k_s\}$.}

\emph{Proof.} Let $k_1, \ldots, k_s$, where $1 \leq s \leq n-2$ are all fixed points of endomorphism $\alpha \in \widehat{\mathcal{E}}_{\mathcal{C}_n}$. Let $\ell \neq k_j$, $j = 1, \ldots, s$ and $\alpha(\ell) = m$. From Proposition 3.1 follows that $\alpha(m) = m$. So, we have  $m = k_j$ for some $j = 1, \ldots, s$. Then $\Im(\alpha) = \{k_1, \ldots, k_s\}$, since $\alpha(k_j) = k_j$ for any $j = 1, \ldots, s$.

Conversely, $\Im(\alpha) = \{k_1, \ldots, k_s\}$, where  $k_1, \ldots, k_s$ are  fixed points of  endomorphism $\alpha$. Now $\alpha(\ell) = k_j$ implies $\alpha(k_j) = k_j$ for arbitrary $j = 1, \ldots, s$, so, using Proposition 3.1 we see that $\alpha$ is an idempotent endomorphism.

\vspace{3mm}

Let $\alpha \in \widehat{\mathcal{E}}_{\mathcal{C}_n}$ have just $n-1$ fixed points. Then $\alpha \neq \mathbf{i}$ and $n - 1 \leq |\Im(\alpha)| < n$. So, $|\Im(\alpha)| = n - 1$ and from Corollary 3.3 follows

\vspace{3mm}

\textbf{Corollary 3.4 } \textsl{Every endomorphism $\alpha \in \widehat{\mathcal{E}}_{\mathcal{C}_n}$ with $n - 1$ fixed points  is an idempotent.}

\vspace{3mm}

Let $\alpha \in \widehat{\mathcal{E}}_{\mathcal{C}_n}$ have just $n-2$ fixed points. Let $\alpha(k) \neq k$ and $\alpha(\ell) \neq \ell$, where $k, \ell \in \mathcal{C}_n$. Assume that $k$ and $\ell$ are not consecutive. Then $\alpha(k) = k- 1$ or $\alpha(k) = k+1$. Similarly $\alpha(\ell) = \ell- 1$ or $\alpha(\ell) = \ell+1$. Since $k-1$, $k+1$, $\ell - 1$ and $\ell + 1$ are fixed points, from Proposition 3.1 follows that $\alpha$ is an idempotent. Assume that $\ell = k + 1$. If
$$\begin{array}{ccc}1. & \alpha(k) = k-1, & \alpha(k+1) = k - 1\\
2. & \alpha(k) = k-1, & \alpha(k+1) = k + 2\\
3. & \alpha(k) = k+2, & \alpha(k+1) = k + 2
\end{array}$$
then easily follows that $\alpha$ is an idempotent. Assume that $\alpha(k) = k + 1$ and $\alpha(k+1) = k+2$. Then from Proposition 3.1 follows that $\alpha$ is not an idempotent endomorphism. But it is clear that $\alpha^2(k) = k+2 = \alpha^3(k)$ and $\alpha^2(k+1) = k + 2 = \alpha^3(k+1)$. Thus we prove

\vspace{3mm}

\textbf{Corollary 3.5 } \textsl{Every endomorphism $\alpha \in \widehat{\mathcal{E}}_{\mathcal{C}_n}$ with $n - 2$ fixed points  is either an idempotent, or $\alpha^3 = \alpha^2$.}

\vspace{3mm}

By similar reasonings we may prove the more general fact

\vspace{3mm}

\textbf{Corollary 3.6 } \textsl{Let $\alpha \in \widehat{\mathcal{E}}_{\mathcal{C}_n}$ have $s$ fixed points $k_1, \ldots, k_s$, where $\ds \left[\frac{n+1}{2}\right] \leq s \leq n-1$ and all other two points $k$ and $\ell$, $k < \ell$, are not consecutive, that is $k\neq k_i$, $\ell \neq k_i$, where $i = 1, \ldots, s$ and $\ell \neq k + 1$. Then $\alpha$   is an idempotent.}

\vspace{3mm}

Let $\alpha \in \widehat{\mathcal{E}}_{\mathcal{C}_n}$ have just two fixed points $k$ and  $\ell$, where $k < \ell$. From Corollary 3.3 follows that
$$\alpha  = \wr\, k, \ldots, k,i_{k+1}, \ldots, i_{\ell-1},\ell, \ldots, \ell\, \wr,$$
where  images $i_s$, $s = k+1, \ldots, \ell -1$, are equal either to  $k$, or to $\ell$.

Thus the endomorphisms of such type are:
$$\alpha_0 = \wr\, k, \ldots, k,  \ldots, k, \ell, \ldots, \ell\, \wr\;; $$

\vspace{-1mm}

{\scriptsize \centerline{$\hphantom{aaaaaaaaaaaaaaaaaaaaaa..}\uparrow \hphantom{aaaaaaaa} \uparrow \hphantom{..aaaaaaaaaaaaaaaa}$ }

\centerline{$\hphantom{.aaaaaaaaaaaaaaaaaaa.a.}k \hphantom{aaaaaaaa.} \ell \hphantom{.aaaaaaaaaaaaaa..}$ } }

\centerline{......................................................$\hphantom{....}$}

$$\alpha_{\ell-k-2} = \wr\, k, \ldots, k, k, \ell  \ldots, \ell, \ldots, \ell\, \wr\;; \hphantom{kk} $$

\vspace{-0.5mm}

{\scriptsize \centerline{$\hphantom{aaaaaaaaaaaaaaaaaaaaaa..}\uparrow \hphantom{aaaaaaaaa,} \uparrow \hphantom{..aaaaaaaaaaaaaaa}$ }

\centerline{$\hphantom{aaaaaaaaaaaaaaaaaaa.}k \hphantom{aaaaaaaaaa.} \ell \hphantom{aaaaaaaaaaaa.}$ } }

$$\alpha_{\ell-k-1} = \wr\, k, \ldots, k, \ell  \ldots, \ell, \ldots, \ell\, \wr. \hphantom{kkkk} $$

\vspace{-1mm}

{\scriptsize \centerline{$\hphantom{aaaaaaaaaaaaaaaaaaaaaa..}\uparrow \hphantom{.aaaaaa} \uparrow \hphantom{..aaaaaaaaaaaaaaaaaa}$ }

\centerline{$\hphantom{aaaaaaaaaaaaaaaaaaa.}k \hphantom{.aaaaaaa} \ell \hphantom{aaaaaaaaaaaaaaa.}$ } }

\vspace{1mm}

For these maps we have
$$\alpha_0 < \ldots < \alpha_{l-k-2} < \alpha_{l-k-1}.$$
Hence $\alpha_i + \alpha_j = \alpha_j + \alpha_i = \alpha_j$, where $j \geq i$. It is easy to see that $\alpha_i\cdot \alpha_j = \alpha_i$ for all $i,j \in \{0, \ldots l-k-2\}$. So, we prove

\vspace{3mm}

\textbf{Proposition 3.7 } \textsl{The set of idempotent endomorphisms with two fixed points $k$ and $\ell$, $k < \ell$, is a semiring of order $\,l - k$ .}

\vspace{3mm}

\textbf{Example 3.8 } The idempotent endomorphisms of semiring à $\widehat{\mathcal{E}}_{\mathcal{C}_7}$ with fixed points 1 and 5 are $$\varphi_1 = \wr\,1\, 1\,1\, 1\, 1\, 5\, 5\, \wr,\; \varphi_2 = \wr\,1\, 1\,1\, 1\, 5\, 5\, 5\, \wr,\; \varphi_3 = \wr\,1\, 1\,1\, 5\, 5\, 5\, 5\, \wr \; \mbox{and}\; \varphi_4 = \wr\,1\, 1\,5\, 5\, 5\, 5\, 5\, \wr.$$
The semiring consisting of these maps has the following addition and multiplication tables
$$
    \begin{array}{c|cccc}
      + & \varphi_1 & \varphi_2 & \varphi_3 & \varphi_4 \\ \hline
      \varphi_1 & \varphi_1 & \varphi_2 & \varphi_3 & \varphi_4 \\
      \varphi_2 & \varphi_2 & \varphi_2 & \varphi_3 & \varphi_4 \\
       \varphi_3 & \varphi_3 & \varphi_3 & \varphi_3 & \varphi_4 \\
        \varphi_4 & \varphi_4 & \varphi_4 & \varphi_4 & \varphi_4 \\
    \end{array}\,,\;\;\;
    \begin{array}{c|cccc}
      + & \varphi_1 & \varphi_2 & \varphi_3 & \varphi_4 \\ \hline
      \varphi_1 & \varphi_1 & \varphi_1 & \varphi_1 & \varphi_1 \\
      \varphi_2 & \varphi_2 & \varphi_2 & \varphi_2 & \varphi_2 \\
       \varphi_3 & \varphi_3 & \varphi_3 & \varphi_3 & \varphi_3 \\
        \varphi_4 & \varphi_4 & \varphi_4 & \varphi_4 & \varphi_4 \\
    \end{array} .
$$

\vspace{4mm}

\textbf{Theorem 3.9 } \textsl{The subset of $\widehat{\mathcal{E}}_{\mathcal{C}_n}$, $n \geq 3$, of all idempotent endomorphisms with $s$ fixed points $k_1, \ldots, k_s$, ${1 \leq s \leq n-1}$, is a semiring of order  $\ds \prod_{m=1}^{s-1} (k_{m+1} - k_{m})$.}

\emph{Proof.} Let $k_1$ be the least fixed point of the idempotent endomorphism $\alpha$. If $k < k_1$, then $\alpha(k) \leq k_1$. Assume that $\alpha(k) = j < k_1$. Then $\alpha(j) = j$ is a contradiction to the minimal choice of $k_1$. Hence $\alpha(k) = k_1$, that is in the first $k_1 + 1$ positions of the ordered  $n$ -- tuple, which represent the endomorphism $\alpha$, occurs only $k_1$.

 Let $k_s$ be the biggest fixed point of the idempotent endomorphism $\alpha$. If $k > k_s$, then $\alpha(k) \geq k_s$. Assume that $\alpha(k) = j > k_s$. Now follows $\alpha(j) = j$ which is a contradiction to the maximal choice of $k_s$. Hence $\alpha(k) = k_s$, that is the last $n - k_s$ positions of the ordered  $n$ -- tuple, which represent the endomorphism $\alpha$, occurs only $k_s$.

 Denote by $p_{k_m}$, where $1 \leq m \leq s$, the number of coordinates in the ordered $n$ -- tuple which represents the endomorphism $\alpha$ equal to $k_m$. For $p_{k_m}$ follows

\vspace{3mm}

\centerline{$k_1 + 1 \leq p_{k_1} \leq k_2$}

\vspace{2mm}

\centerline{$\hphantom{aaaaaaaa.}1  \leq p_{k_2} \leq k_3 - k_2$}

\vspace{2mm}

\centerline{$\hphantom{aaaaaaaa}$..........................................}

\vspace{2mm}

\centerline{$\hphantom{aaaaaaaaaaaa.}1  \leq p_{k_m} \leq k_{m+1} - k_m$  }

\vspace{2mm}

\centerline{$\hphantom{aaaaaaaaa}$.............................................}

\vspace{2mm}

\centerline{$\hphantom{aaaaaaaaa}n - k_{s}  \leq p_{k_{s}}\, \leq\, n - k_{s- 1} - 1$}

\vspace{3mm}

So every idempotent endomorphism with $s$ fixed points $k_1, \ldots, k_s$ have the form
$$\alpha = \wr \, \underbrace{k_1, \ldots, k_1}, \underbrace{k_2, \ldots, k_2}, \ldots, \underbrace{k_s, \ldots, k_s}\, \wr $$

 \vspace{-3mm}

\centerline{$\hphantom{aaa..aa}p_{k_1}$ $\hphantom{aaaa.a}p_{k_2}$ $\hphantom{aaaaaaaaa}p_{k_s}$ $\hphantom{.}$}
\vspace{2mm}

Let another endomorphism $\beta$ has the same form

$$\beta = \wr \, \underbrace{k_1, \ldots, k_1}, \underbrace{k_2, \ldots, k_2}, \ldots, \underbrace{k_s, \ldots, k_s}\, \wr $$

 \vspace{-1mm}

\centerline{$\hphantom{aaa..aa}p'_{k_1}$ $\hphantom{aaaa.a}p'_{k_2}$ $\hphantom{aaaaaaaaa}p'_{k_s}$ $\hphantom{.}$}

\vspace{2mm}

Since $(\alpha + \beta)(k_i) = \alpha(k_i) + \beta(k_i) = k_i + k_i = k_i$, where $i = 1, \ldots, s$, follows that $\alpha + \beta$ has the form of $\alpha$ and $\beta$.

For these $\alpha$ and $\beta$ follows $\alpha\cdot \beta = \alpha$. (So, the multiplication of the idempotent endomorphisms with fixed points $k_1, \ldots, k_s$ is like the multiplication of the constant endomorphisms -- every idempotent is a right identity.)

Thus  we prove that the set of all idempotent endomorphisms with $s$ fixed points $k_1, \ldots, k_s$, ${1 \leq s \leq n-1}$, is a semiring.

From the inequalities of numbers $p_{k_m}$, where $1 \leq m \leq s $, follows that all  possibilities to have $k_m$ in $n$ -- tuple are $k_{m+1} - k_m$.

The order of this semiring is equal to the product of all possibilities for the fixed points ${k_m}$, where $1 \leq m \leq s$. This number is equal to $\ds \prod_{m=1}^{s-1} (k_{m+1} - k_{m})$ and the proof is completed.

\vspace{3mm}

The semiring of the idempotent endomorphisms of $\widehat{\mathcal{E}}_{\mathcal{C}_n}$ with $s$ fixed points  $k_1, \ldots, k_s$ is denoted by $\mathcal{ID}(k_1, \ldots, k_s)$. Since the semiring of all endomorphisms with fixed points $k_1, \ldots, k_s$ is $\ds \bigcap_{r = 1}^s \mathcal{E}^{(k_r)}_{\mathcal{C}_n}$, it follows that $\mathcal{ID}(k_1, \ldots, k_s)$ is a subsemiring of $\ds \bigcap_{r = 1}^s \mathcal{E}^{(k_r)}_{\mathcal{C}_n}$.

\vspace{3mm}

\textbf{Remark 3.10 } a. From Theorem 3.9  semiring  $\mathcal{ID}(0, \ldots, k - 1, k + 1, \ldots, n - 1)$, has two elements: the idempotent endomorphisms $k^+ = \wr \, 0, \ldots, k-1, k+1, k+1, \ldots n-1\, \wr$ and $k^- = \wr \, 0, \ldots, k-1, k-1, k+1, \ldots n-1\, \wr$. This semiring has the following addition and multiplication tables
$$\begin{array}{l|ll}
    + & k^- &  k^+ \\ \hline
    k^- & k^- &  k^+ \\
    k^+ & k^+ &  k^+
  \end{array}\, , \;\;\;
\begin{array}{l|ll}
    \cdot & k^- &  k^+ \\ \hline
    k^- & k^- &  k^- \\
    k^+ & k^+ &  k^+
  \end{array}.
$$

b. The product of endomorphisms from different semirings $\mathcal{ID}(k_1, \ldots, k_s)$ is not, in general, an idempotent. Indeed, in semiring $\widehat{\mathcal{E}}_{\mathcal{C}_8}$ for
$$\alpha = \wr \, 1\, 1\, 5\, 5\, 5\, 5\, 5\, 5\, \wr \in \mathcal{ID}(1,5) \;\; \mbox{and} \;\; g = \wr \, 2\, 2\, 2\, 5\, 5\, 5\, 5\, 5\, \wr \in \mathcal{ID}(2,5)$$
follows $\alpha\cdot \beta = \wr \, 2\, 2\, 5\, 5\, 5\, 5\, 5\, 5\, \wr$, but this endomorphism is not an idempotent.

c. The sum of two idempotent endomorphisms such that the fixed points of the first one  are part of the fixed points of another endomorphism, is not an idempotent. For instance\break
 $\wr\,0\,0\,3\,3\,\wr + \wr\,0\,2\,2\,3\,\wr = \wr\,0\,2\,3\,3\,\wr$.

\vspace{3mm}

\textbf{Corollary 3.11 } \textsl{Let $n \geq 3$ and $k_1, \ldots, k_s \in \mathcal{C}_n$, where $s = 1, \ldots, n-2$. Semiring $\mathcal{ID}(k_1, \ldots, k_s)$ is an ideal of  $\ds \bigcap_{r = 1}^s \mathcal{E}^{(k_r)}_{\mathcal{C}_n}$.}

\emph{Proof.} Using Theorem 3.9 it will be enough to show that $\mathcal{ID}(k_1, \ldots, k_s)$ is closed under left and right multiplications by elements of
 $\ds \bigcap_{r = 1}^s \mathcal{E}^{(k_r)}_{\mathcal{C}_n}$.

 Let $$\alpha = \wr \, \underbrace{k_1, \ldots, k_1}, \underbrace{k_2, \ldots, k_2}, \ldots, \underbrace{k_s, \ldots, k_s}\, \wr \in \mathcal{ID}(k_1, \ldots, k_s)\; \mbox{and}\;$$

 \vspace{-0.5mm}

\centerline{$\hphantom{aaaa..}p_{k_1}$ $\hphantom{aaaa.a}p_{k_2}$ $\hphantom{aaaaaaaaa}p_{k_s}$ $\hphantom{.aaaaaaaaaaaaaaaaa}$}

 \vspace{-2mm}

$$\beta = \wr \, i_0, i_1, \ldots, i_{k_1 -1},k_1, i_{k_1 +1}, \ldots, i_{k_2 -1}, k_2, \ldots, k_s\, \wr \, \in \; \bigcap_{r = 1}^s \mathcal{E}^{(k_r)}_{\mathcal{C}_n}.\hphantom{aa}$$
Hence $\alpha \cdot \beta = \alpha$.
Calculate
$$\beta\cdot \alpha = \wr \, i_0, i_1, \ldots, i_{k_1 -1},k_1, i_{k_1 +1}, \ldots, i_{k_2 -1}, k_2, \ldots, k_s\, \wr\cdot \wr \, {k_1, \ldots, k_1}, {k_2, \ldots, k_2}, \ldots, {k_s, \ldots, k_s}\, \wr.$$

Since $i_s \leq k_1$, where $s = 0, \ldots, k_1 - 1$, then $\alpha(i_s) \leq \alpha(k_1) = k_1$ and from the proof of Theorem 3.9 it  follows that $\alpha(i_s) = k_1$.

Since $k_1 \leq i_{k_1 + 1} \leq \cdots \leq i_{k_2 - 1} \leq k_2$, it follows
$$k_1 = \alpha(k_1) \leq \alpha\left(i_{k_1 + 1}\right) \leq \cdots \leq \alpha\left(i_{k_2 - 1}\right) \leq \alpha(k_2) = k_2.$$

But endomorphism $\alpha$ maps all elements between $k_1$ and $k_2$ either in $k_1$, or in $k_2$. Then $\alpha\left(i_{k_1 + 1}\right), \ldots, \alpha\left(i_{k_2 - 1}\right)$ are either $k_1$, or $k_2$. We use such arguments for the next elements of $\Im(\beta)$ which are between other fixed points.
Thus
$$\beta\cdot \alpha = \gamma = \wr \, \underbrace{k_1, \ldots, k_1}, \underbrace{k_2, \ldots, k_2}, \ldots, \underbrace{k_s, \ldots, k_s}\, \wr \in \mathcal{ID}(k_1, \ldots, k_s).$$

 \vspace{-1.5mm}

\centerline{$\hphantom{aaaaaaa..}p'_{k_1}$ $\hphantom{aaaa.a}p'_{k_2}$ $\hphantom{aaaaaaaaa}p'_{k_s}$ $\hphantom{.aaaaaaaaaa.}$}

\vspace{7mm}

\centerline{{\large 4. \hspace{0.5mm} Roots of idempotent endomorphisms}}

\vspace{3mm}

Let $(S,\cdot)$ be a finite semigroup. It is well known, see [6], that for any $x \in S$ there is a positive integer $k = k(x)$ such that $a^k$ is an idempotent element of $S$.

Now we consider the following relation:
For any $x, y \in S$ define $$x \sim y \iff \exists\, k, m \in \mathbf{N},\; x^k = y^m = e,$$
where $e$ is an idempotent element of $S$.

Obviously, the relation $\sim$ is reflexive and symmetric.

Let $x \sim y$ and $y \sim z$. Then  $x^k = y^m = e_1$ and $y^r = z^s = e_2$ where $e_1$ and $e_2$ are idempotents. Now it follows $\left(y^m\right)^r = e_1^r = e_1$ and $\left(y^r\right)^m = e_2^m = e_2$, thus $e_1 = e_2$. Hence $x^k = z^s = e_1$, that is $x \sim z$. So, we prove that $\sim$ is an equivalence relation on $S$.

Note that two different idempotents belong to different equivalence classes modulo $\sim$. If $e$ is an idempotent the elements of the equivalence class containing $e$ are called \emph{roots of the idempotent} $e$.

\vspace{1mm}

The following natural questions arise:

\emph{
1. Are there any finite semigroups where the above relation is a congruence?}

\emph{
2. Are there any finite semirings where the above relation is a congruence?}

\emph{
3. Are there equivalence classes which are semigroups?}

\vspace{3mm}

To answer the third question we consider semigroup $\left(\widehat{\mathcal{E}}_{\mathcal{C}_n},\cdot\right)$ and the defined above equivalence relation.

Let $\alpha \in \widehat{\mathcal{E}}_{\mathcal{C}_n}$. Element $\mathbf{j}  \in \mathcal{C}_n$ is called jump point of $\alpha$ if $\mathbf{j} \neq 0$ and one of the following conditions hold:

 1. $\alpha(\mathbf{j} - 1) \leq \mathbf{j} - 1$ and $\alpha(\mathbf{j}) > \mathbf{j}$,

 2. $\alpha(\mathbf{j} - 1) < \mathbf{j} - 1$ and $\alpha(\mathbf{j}) \geq \mathbf{j}$.

There are endomorphisms without jump points, namely,  identity $\mathbf{i}$ and  constant endomorphisms $\kappa_k = \wr\,k \,k\, \ldots\, k\,\wr$, $k = 0, 1, \ldots, n-1$. For $k, \ell  \in \mathcal{C}_n$, $k < \ell$  endomorphism $\ds \alpha_j(i) = \left\{\begin{array}{ll} k, & i \leq j - 1\\ \ell, & i \geq j \end{array} \right.$ has no jump points if $j > \ell$. Endomorphisms $k^+$ and $k^-$, see Remark 3.10, has jump points: $k$ and $k+1$ are jump points of $k^+$ and $k^-$, respectively.

\vspace{3mm}

\textbf{Theorem 4.1 } \textsl{Let  $\alpha \in \widehat{\mathcal{E}}_{\mathcal{C}_n}$, $n \geq 3$, is an endomorphism  with $s$ fixed points $k_1, \ldots, k_s$, ${1 \leq s \leq n-2}$. Let for some $i$, $i = 1, \ldots, s -1$,  fixed points $k_i$ and $k_{i+1}$ be not consecutive, i.e. $k_{i+1} \neq k_i + 1$. Then there is a unique jump point $\mathbf{j}_{\,i}$ of $\alpha$ such that $k_i+1 \leq \mathbf{j}_{\,i} \leq k_{i+1}$. }

\vspace{1mm}

\emph{Proof.} For $k_i + 1$ there are two possibilities:

\vspace{1mm}

1. $\alpha(k_i + 1) > k_i + 1$, then $\mathbf{j}_{\,i} = k_i + 1$ is the searched  jump point of $\alpha$, or

2. $\alpha(k_i + 1) < k_i + 1$.

Let $\alpha(k_i + m) < k_i + m$ for any $m = 1, 2, \ldots, \ell - 1$. If $\alpha(k_i + \ell) \geq k_i + \ell$, then $\mathbf{j}_{\,i} = k_i + \ell$ is the searched  jump point of $\alpha$.

Assume that $\alpha(k_i + m) < k_i + m$ for all $m = 1, 2, \ldots, k_{i+1} - k_i - 1$. Now $\alpha(k_{i+1} - 1) < k_{i+1} - 1$. But $\alpha(k_{i+1}) = k_{i+1}$, so, $\mathbf{j}_{\,i} = k_{i+1}$ is the searched  jump point of $\alpha$.

Suppose that $\mathbf{j}$ and $\mathbf{j}^{\,\circ}$ are two jump points of $\alpha$ such that $k_i+1 \leq \mathbf{j} \leq k_{i+1}$ and\break $k_i+1 \leq \mathbf{j}^{\,\circ} \leq k_{i+1}$. Let $\mathbf{j} < \mathbf{j}^{\,\circ}$. Now $\alpha(\mathbf{j}) > \mathbf{j}$. Let $m \in \mathcal{C}_n$, $\mathbf{j} \leq m < \mathbf{j}^{\,\circ}$ is the maximal element such that $\alpha(m) > m$. Then $\alpha(m+ 1) < m + 1$ and it follows $\alpha(m) \geq m + 1 > \alpha(m + 1)$ which is a contradiction. By the same arguments we show that $\mathbf{j}^{\,\circ} < \mathbf{j}$ is impossible. So $\mathbf{j} = \mathbf{j}^{\,\circ}$ is the unique jump point of $\alpha$ such that $k_i+1 \leq \mathbf{j} \leq k_{i+1}$.

\vspace{3mm}

Let $k_1$ be the least fixed point of  endomorphism $\alpha$ and $i \in \mathcal{C}_n$, $i < k_1$. If suppose that $\alpha(i) < i$, using $\alpha(0) > 0$, it follows that for some maximal element $m$, where $0 \leq m < i$, we have $\alpha(m) > m$ and then $\alpha(m+1) < m+1$, which is impossible. So, for every $i < k_1$ it follows $\alpha(i) > i$. Similarly for every $i > k_s$, where $k_s$ is the biggest fixed point of $\alpha$, it follows that $\alpha(i) < i$. So, we prove

\vspace{3mm}

\textbf{Corollary 4.2 } \textsl{Every endomorphism  $\alpha \in \widehat{\mathcal{E}}_{\mathcal{C}_n}$, $n \geq 3$,   with just $\ell$ fixed points $k_1, \ldots, k_\ell$, which are not consecutive, i.e. $k_{i+1} \neq k_i + 1$ for $i = 1, \ldots, \ell - 1$, has just $\ell - 1$ jump points $\mathbf{j}_{\,i}$ such that $k_i+1 \leq \mathbf{j}_{\,i} \leq k_{i+1}$.}

\vspace{3mm}

\textbf{Remark 4.3 } From Theorem 3.9 it follows that the number of idempotents with fixed points $k_1, \ldots, k_s$ is a semiring of order  $\ds \prod_{i=1}^{s-1} (k_{i+1} - k_{i})$.
But if two fixed points $k_{i}$ and $k_{i+1}$ are consecutive, the difference $k_{i+1} - k_i$ is equal to 1 and they do not appear in the product. So, the order of semiring is equal to $\ds \prod_{i=1}^{\ell-1} (k_{i+1} - k_{i})$, where $k_1, \ldots, k_\ell$ are (after suitable renumbering) not consecutive fixed points. But for any $i$ the difference $k_{i+1} - k_i$ is just the number of possible jump points $\mathbf{j}_{\,i,t} = k_i + t$, where $t = 1, \ldots, k_{i+1} - k_i$. So, for given fixed points $k_1, \ldots, k_s$ the order of the semiring of idempotents is equal to the number of all $\ell - 1$ -- tuples $(\mathbf{j}_{\,1,t_1}, \mathbf{j}_{\,2,t_2}, \ldots, \mathbf{j}_{\,\ell - 1,t_{\ell-1}})$, where  $\mathbf{j}_{\,i,t_i} = k_{i} + t_i$ for $t_i = 1, \ldots, k_{i+1} - k_i$ and  $i = 1, \ldots, \ell - 1$. Hence, any $\ell - 1$ -- tuple $(\mathbf{j}_{\,1,t_1},  \ldots, \mathbf{j}_{\,\ell - 1,t_{\ell-1}})$ define just one idempotent from the considered semiring.

\vspace{3mm}

To describe precisely all the fixed points of an arbitrary endomorphism $\alpha \in \widehat{\mathcal{E}}_{\mathcal{C}_n}$ we shall use new indices. Let first fixed point of $\alpha$ is $k_{1,1}$ and some fixed points after $k_{1,1}$ be  consecutive, i.e. $k_{1,2}, \ldots, k_{1,m_1}$ are fixed points such that
 $k_{1,\,i+1} - k_{1,\,i} = 1$ where $i = 1, \ldots, m_1 - 1$. Let the next fixed point be $k_{2,1}$ such that $k_{2,1} - k_{1,m_1} > 1$. So we construct the first pair of two fixed points which are not  consecutive. Let the following fixed points are $k_{2,2}, \ldots, k_{2,m_2}$ such that  $k_{2,\,i+1} - k_{2,\,i} = 1$ where $i = 1, \ldots, m_2 - 1$. The next fixed point is $k_{3,1}$ and $k_{3,1} - k_{2,m_2} > 1$. Let the last pair of two not consecutive fixed points be $k_{\ell-1,\,m_{\ell-1}}$ and $k_{\ell,\,1}$. Then the last fixed points are $k_{\ell,\,2}, \ldots k_{\ell,\,m_\ell}$ such that  $k_{\ell,\,i+1} - k_{\ell,\,i} = 1$ where $i = 1, \ldots, m_\ell - 1$. So, we construct a partition of a set of fixed points of $\alpha$ such that we may distinguish the fixed points which are not consecutive.

Let $\mathbf{j}_{\,i,t_i}$ be the jump points of $\alpha$ such that $\mathbf{j}_{\,i,t_i} = k_{i,\,m_i} + t_i$, where $t_i = 1, \ldots, k_{i+1,\,1} - k_{i,\,m_i}$ and $i = 1, \ldots, \ell - 1$.

An endomorphism $\alpha$ with fixed points $k_{1,1}, \ldots, k_{\ell,\,m_\ell}$ and jump points $\mathbf{j}_{\,i,t_i}$ from the previous definitions is called endomorphism of type $$[k_{1,1}, \ldots, k_{1,m_1},\mathbf{j}_{\,1,t_1},k_{2,1}, \ldots, k_{\ell-1,\,m_{\ell-1}},\mathbf{j}_{\,\ell-1,t_{\ell-1}},k_{\ell,1}, \ldots, k_{\ell,\,m_\ell}]. \eqno{(1)}$$

Let us consider  endomorphism $\varepsilon$ of this type such that:

$\bullet$ $\varepsilon(x) = k_{1,1}$ for any $0 \leq x \leq k_{1,1}$,

$\bullet$ $\varepsilon(x) = k_{i,\,m_i}$ for any $k_{i,\,m_i} \leq x \leq \mathbf{j}_{\,i,t_i} - 1$, where $i = 1, \ldots, \ell - 1$,

$\bullet$ $\varepsilon(x) = k_{i+1,1}$ for any $\mathbf{j}_{\,i,t_i} \leq x \leq k_{i+1,\,1}$, where $i = 1, \ldots, \ell - 1$,

$\bullet$ $\varepsilon(x) = k_{\ell,\,m_\ell}$ for any $k_{\ell,\,m_\ell} \leq x \leq n - 1$.

Now it is  easy to show that this endomorphism is an idempotent.

Let $\bar{\varepsilon}$ is another idempotent of the same type. Then using the reasonings just before Corollary 5.2 follows that

$\bullet$ $\bar{\varepsilon}(x) = k_{1,1}$ for any $0 \leq x \leq k_{1,1}$,

$\bullet$ $\bar{\varepsilon}(x) = k_{\ell,\,m_\ell}$ for any $k_{\ell,\,m_\ell} \leq x \leq n - 1$.

Since $\bar{\varepsilon}$ is an idempotent we conclude that for some $x$, where $k_{i,\,m_i} \leq x \leq k_{i+1,\,1}$, it follows either $\bar{\varepsilon}(x) = k_{i,\,m_i}$, or $\bar{\varepsilon}(x) = k_{i + 1,\,1}$. But using that $\bar{\varepsilon}$ is of type (1) it follows that $\bar{\varepsilon}(x) = \varepsilon(x)$ for all $x$, where $k_{i,\,m_i} \leq x \leq k_{i+1,\,1}$ and  $i = 1, \ldots, \ell -1$. Hence, there is only one idempotent of a given type (1). This endomorphism is
$$\varepsilon = \wr\,k_{1,1}, \ldots, k_{1,1},k_{1,2}, \ldots, k_{1,\,m_1}, \ldots, k_{1,\,m_1},k_{2,1}, \ldots, k_{\ell-1,\,m_{\ell-1}},k_{\ell,\,1}, \ldots, k_{\ell,\,m_\ell}, \ldots, k_{\ell,\,m_\ell}\wr. \eqno{(2)}$$

\vspace{-2mm}

{\small $\;$ \hspace{17.8mm} $\uparrow\;\;$\hspace{17.9mm} $\uparrow\;\;$ \hspace{22.9mm} $\uparrow\;\;$ \hspace{26.7mm} $\uparrow\;\;$\hspace{10mm} $\uparrow\;\;$ }

{\small $\;$ \hspace{17.8mm} $k_{1,1}$\hspace{15.9mm} $k_{1, m_1}$ \hspace{18.5mm} $\mathbf{j}_{\,1,t_1}$ \hspace{23.7mm} $\mathbf{j}_{\,\ell-1,t_{\ell-1}}$ \hspace{-0.5mm} $k_{\ell, m_\ell}$}

\vspace{5mm}

\textbf{Lemma 4.4 } \textsl{Let $\varepsilon \in \widehat{\mathcal{E}}_{\mathcal{C}_n}$ be an idempotent endomorphism. If $\alpha \in \widehat{\mathcal{E}}_{\mathcal{C}_n}$ is a root of $\varepsilon$, then the two endomorphisms $\alpha$ and $\varepsilon$ have the same fixed points.
}

\emph{Proof.} Let $k$ be a fixed point of $\varepsilon$. If $\alpha(k) < k$, then $\alpha^2(k) \leq \alpha(k) < k$ and by the same arguments $\alpha^m(k) < k$ for arbitrary natural number $m$, but this is a contradiction with $\alpha^m = \varepsilon$ for some $m \in \mathbf{N}$. Analogously, assume that $\alpha(k) > k$, we obtain a contradiction to equality $\alpha^m = \varepsilon$. Hence $\alpha(k) = k$.

Let $k$ be a fixed point of $\alpha$. Since it is a fixed point of $\alpha^m$ for arbitrary natural $m$ it follows that it is a fixed point of $\varepsilon$.

\vspace{3mm}

\textbf{Lemma 4.5 } \textsl{Let $\varepsilon \in \widehat{\mathcal{E}}_{\mathcal{C}_n}$ be an idempotent endomorphism. If $\alpha \in \widehat{\mathcal{E}}_{\mathcal{C}_n}$ is a root of $\varepsilon$, then the two endomorphisms $\alpha$ and $\varepsilon$ have the same jump points.}

\emph{Proof.} Let $\mathbf{j}$ be a jump point of $\varepsilon$ and $\varepsilon(\mathbf{j} - 1) \leq \mathbf{j} - 1$, $\varepsilon(\mathbf{j}) > \mathbf{j}$. If we assume that $\alpha(\mathbf{j} - 1) > \mathbf{j} - 1$, then $\alpha^m(\mathbf{j} - 1) > \mathbf{j} - 1$ for arbitrary natural $m$ which is  a contradiction with $\alpha^m = \varepsilon$ for some $m \in \mathbf{N}$.

 Analogously, if we assume $\alpha(\mathbf{j}) \leq \mathbf{j}$ it follows $\alpha^m(\mathbf{j}) \leq \mathbf{j}$ for arbitrary natural $m$ which is also a contradiction. Hence $\mathbf{j}$ is a jump point of $\alpha$ such that $\alpha(\mathbf{j} - 1) \leq \mathbf{j} - 1$ and $\alpha(\mathbf{j}) > \mathbf{j}$. By similar reasonings it follows that if $\mathbf{j}$ is a jump point of $\varepsilon$ such that $\varepsilon(\mathbf{j} - 1) < \mathbf{j} - 1$ and $\varepsilon(\mathbf{j}) \geq \mathbf{j}$ for every endomorphism $\alpha$, which is root of $\varepsilon$, we have that $\mathbf{j}$ is a jump point of $\alpha$, $\alpha(\mathbf{j} - 1) < \mathbf{j} - 1$ and $\alpha(\mathbf{j}) \geq \mathbf{j}$.

Let $\mathbf{j}$ be a jump point of $\alpha$ such that $\alpha(\mathbf{j} - 1) \leq \mathbf{j} - 1$ and $\alpha(\mathbf{j}) > \mathbf{j}$. Then for arbitrary natural $m$ follows $\alpha^m(\mathbf{j} - 1) \leq \mathbf{j} - 1$ and $\alpha^m(\mathbf{j}) > \mathbf{j}$. So, $\mathbf{j}$ is a jump point of $\varepsilon$,  $\varepsilon(\mathbf{j} - 1) \leq \mathbf{j} - 1$ and $\varepsilon(\mathbf{j}) > \mathbf{j}$. Analogously $\alpha(\mathbf{j} - 1) < \mathbf{j} - 1$ and $\alpha(\mathbf{j}) \geq \mathbf{j}$ implies $\alpha^m(\mathbf{j} - 1) < \mathbf{j} - 1$ and $\alpha^m(\mathbf{j}) \geq \mathbf{j}$ for every natural number $m$ and so $\mathbf{j}$ is a jump point of $\varepsilon$,  $\varepsilon(\mathbf{j} - 1) < \mathbf{j} - 1$ and $\varepsilon(\mathbf{j}) \geq \mathbf{j}$.

\vspace{3mm}

Immediately from lemmas 4.4 and 4.5 follows

\vspace{3mm}

\textbf{Proposition 4.6 } \textsl{All the endomorphisms of one equivalence class modulo $\sim$ have the same type.}

\vspace{3mm}

Let the idempotent endomorphism $\varepsilon$ from (2) be an element of  equivalence class $E$. Then $E$ is called \emph{equivalence class of type} (1).

\vspace{3mm}

\textbf{Lemma 4.7 } \textsl{Let $E$ be an equivalence class of type (1). Then any $\alpha \in E$ satisfies the following conditions:}

 \textbf{A.} $\alpha(x) > x$, where $0 \leq x < k_{1,1}$, or $\mathbf{j}_{\,1,t_1} \leq x < k_{2,1}$, or $\mathbf{j}_{\,2,t_2} \leq x < k_{3,1}$,  $\ldots$, or $\mathbf{j}_{\,\ell-1,t_{\ell-1}} \leq x < k_{\ell,\,1}$.

\textbf{B.} $\alpha(x) < x$, where $k_{1,\,m_1} < x < \mathbf{j}_{\,1,t_1}$, or $k_{2,\,m_2} < x < \mathbf{j}_{\,2,t_2}$, $\ldots$, or
 $k_{\ell-1,\,m_{\ell-1}} < x < \mathbf{j}_{\,\ell-1,t_{\ell-1}}$, or $k_{\ell,\,m_\ell} < x \leq n-1$.

\vspace{1mm}

\emph{Proof.} Assume that for some $x$ such that
$$0 \leq x < k_{1,1}, \; \mbox{or}\; \mathbf{j}_{\,1,t_1} \leq x < k_{2,1}, \; \mbox{or}\; \mathbf{j}_{\,2,t_2} \leq x < k_{3,1}, \ldots, \; \mbox{or}\; \mathbf{j}_{\,\ell-1,t_{\ell-1}} \leq x < k_{\ell,\,1}$$
it follows $\alpha(x) \leq x$. Then $\alpha^2(x) \leq \alpha(x) \leq x$ and by induction $\alpha^m(x) \leq x$ for arbitrary natural $m$. It is a contradiction to the assumption that $\alpha^m = \varepsilon$ for some natural $m$, where $\varepsilon$ has the form (2). Thus we prove that $\alpha(x) > x$.

Analogously, assume that for some $x$ such that
$$k_{1,\,m_1} < x < \mathbf{j}_{\,1,t_1}, \; \mbox{or}\; k_{2,\,m_2} < x < \mathbf{j}_{\,2,t_2}, \ldots, \; \mbox{or}\; k_{\ell-1,\,m_{\ell-1}} < x < \mathbf{j}_{\,\ell-1,t_{\ell-1}},  \; \mbox{or}\; k_{\ell,\,m_\ell} < x \leq n-1$$
 we have $\alpha(x) \geq x$. Then  by induction $\alpha^m(x) \geq x$ for arbitrary natural $m$. This contradicts to the  assumption that $\alpha^m = \varepsilon$ for some natural $m$, where $\varepsilon$ has the form (2). Thus we prove that $\alpha(x) < x$.

\vspace{3mm}

Note that if two endomorphisms have the same fixed points but different jump points their product may posses a new fixed point. For instance  endomorphisms $\alpha =  \wr\, 2,2,2,4,5,5,5\,\wr \sim \varepsilon_1 = \wr\, 2,2,2,5,5,5,5\,\wr$ and $\beta = \wr\, 2,2,2,2,3,5,5\,\wr \sim \varepsilon_2 = \wr\, 2,2,2,2,2,5,5\,\wr$ have the fixed points 2 and 5 but $\alpha\cdot \beta = \wr\, 2,2,2,3,5,5,5\,\wr$ and $\beta\cdot\alpha = \wr\, 2,2,2,2,4,5,5\,\wr$ have the new fixed points 3 and 4, respectively.

\vspace{3mm}

Let us remind that the jump points, which we use in (1) are defined by equality $\mathbf{j}_{\,i,t_i} = k_{i,\,m_i} + t_i$. So it follows $\mathbf{j}_{\,i,t_i} - k_{i,\,m_i} = t_i$. Now we denote $k_{i+1,\,1} - \mathbf{j}_{\,i,t_i} = s_i$. The indices $t_i$ and $s_i$ will be used in the next theorem.

\vspace{4mm}

\textbf{Theorem 4.8 } \textsl{Every equivalence class  modulo $\sim$ of type (1) is a subsemiring of $\widehat{\mathcal{E}}_{\mathcal{C}_n}$, $n \geq 2$. The order of this semiring is
$$C_{k_{1,1}} \prod_{i =1}^{\ell - 1} \left(C_{t_i} C_{s_i}\right)C_{n-1-k_{\ell,\,m_\ell}},$$
where $C_p$ is  $p$ -- th Catalan number.}

\vspace{1mm}

\emph{Proof.} Let $E$ be an equivalence class of type (1) and $\alpha, \beta \in E$.
Using Lemma 4.7 it follows:

\textbf{A.} $\alpha(x) > x$ and $\beta(x) > x$, where $0 \leq x < k_{1,1}$, or $\mathbf{j}_{\,1,t_1} \leq x < k_{2,1}$, or $\mathbf{j}_{\,2,t_2} \leq x < k_{3,1}$,  $\ldots$, or\break $\mathbf{j}_{\,\ell-1,t_{\ell-1}} \leq x < k_{\ell,\,1}$. Then $(\alpha+\beta)(x) = \alpha(x)\vee \beta(x) > x\vee x = x$ and $(\alpha\cdot \beta)(x) = \beta(\alpha(x)) \geq \beta(x) > x$.

\textbf{B.} $\alpha(x) < x$ and $\beta(x) < x$, where $k_{1,\,m_1} < x < \mathbf{j}_{\,1,t_1}$, or $k_{2,\,m_2} < x < \mathbf{j}_{\,2,t_2}$, $\ldots$, or
 $k_{\ell-1,\,m_{\ell-1}} < x < \mathbf{j}_{\,\ell-1,t_{\ell-1}}$, or $k_{\ell,\,m_\ell} < x \leq n-1$.  Then $(\alpha + \beta)(x) = \alpha(x)\vee \beta(x) < x\vee x= x$ and $(\alpha\cdot \beta)(x) = \beta(\alpha(x)) \leq \beta(x) < x$.

Thus we prove that equivalence class $E$ is a subsemiring of $\widehat{\mathcal{E}}_{\mathcal{C}_n}$.

\vspace{3mm}

For the second part of the proof of the theorem we needed  the two combinatorial lemmas.

\vspace{3mm}

\textbf{Lemma 4.9 } \textsl{ The number of the ordered $p$ -- tuples $(i_0, \ldots, i_{p-1})$, where
}

\textsl{1. $i_r \in \{0, \ldots, p-1\}$ for $r = 0, \ldots, p-1$,}

\textsl{2. $i_r \leq i_{r+1}$ for $r = 0, \ldots, p-1$,}

\textsl{3. $i_r < r$ for $r = 1, \ldots, p-1$}

\noindent \textsl{is  $p$ -- th Catalan number $\ds C_p = \frac{1}{p}\binom{2p-1}{p-1}$.}

\vspace{2mm}

The proof of this lemma is a part of the proof of Proposition 3.4 of [10]. This result and many other applications of Catalan numbers can be found in [8].

\vspace{3mm}

\textbf{Lemma 4.10 } \textsl{ The number of the ordered $p$ -- tuples $(i_0, \ldots, i_{p-1})$, where
}

\textsl{1. $i_r \in \{1, \ldots, p\}$ for $r = 0, \ldots, p-1$, $i_p = p$,}

\textsl{2. $i_r \leq i_{r+1}$ for $r = 0, \ldots, p-1$,}

\textsl{3. $i_r > r$ for $r = 0, \ldots, p-1$}

\noindent \textsl{is  $p$ -- th Catalan number $\ds C_p = \frac{1}{p}\binom{2p-1}{p-1}$.}

\vspace{2mm}

\emph{Proof of Lemma 4.10.} Let $(i_0, \ldots, i_{p-1})$ is the ordered $p$ -- tuple.  Since $i_r > r$ for any $r = 0, \ldots, p-1$ follows that $p = i_{p} \geq i_{p-1} > p - 1$, i.e. $i_{p-1} = p$. Instead of the $p$ -- tuple $(i_0, \ldots, i_{p-1})$ we consider the $p$ -- tuple $(j_0, \ldots, j_{p-1})$, where $j_m = p - i_{p-m-1}$, $m = 0, \ldots p-1$. Then $j_m < m$ and $j_m \leq j_{m+1}$. Using Lemma 4.9 follows that the number of all $k$ -- tuples of this kind is $C_k$.

\vspace{3mm}

Now we continue the proof of the theorem.

For all the intervals considered in  part \textbf{A.} (from the begining of this proof) we apply Lemma 4.10. Then, for all the intervals considered in the part \textbf{B.} we apply Lemma 4.9. Hence we find that  the number of the endomorphisms in the equivalence class  modulo $\sim$ of type (1) is
$\ds C_{k_{1,1}} \prod_{i =1}^{\ell - 1}\left(C_{t_i} C_{s_i}\right)C_{n-1-k_{\ell,\,m_\ell}}$.

\vspace{3mm}

\textbf{Remark 4.11 } The  equivalence relation $\sim$ is not a congruence on $\widehat{\mathcal{E}}_{\mathcal{C}_n}$. For $n = 8$ we consider the endomorphisms
$$\begin{array}{l} \alpha =  \wr\, 2,2,2,2,2,6,7,7\,\wr \sim \varepsilon_1 = \wr\, 2,2,2,2,2,7,7,7\,\wr,\\
\beta =  \wr\, 2,2,2,2,3,7,7,7\,\wr \sim \varepsilon_1 = \wr\, 2,2,2,2,2,7,7,7\,\wr,\\
\gamma =  \wr\, 1,1,1,3,3,4,5,6\,\wr \sim \varepsilon_2 = \wr\, 1,1,1,3,3,3,3,3\,\wr.
\end{array}
$$
Now we compute $\alpha\cdot\gamma =  \wr\, 1,1,1,1,1,5,6,6\,\wr$ and $\beta\cdot\gamma =  \wr\, 1,1,1,1,3,6,6,6\,\wr$. So, it follows that $\alpha \sim \beta$, but $\alpha\cdot\gamma \nsim \beta\cdot\gamma$.

\vspace{7mm}

\centerline{{\large 5. \hspace{0.5mm} The crucial role of the jump points}}

\vspace{4mm}

Here we consider idempotent endomorphisms with arbitrary fixed points but assume that each endomorphism have  $\mathbf{j}_1, \ldots, \mathbf{j}_r$ for jump points.

 Two jump points $\mathbf{j}_s$ and $\mathbf{j}_{s+1}$  of the idempotent $\ep$ are called consecutive if $\mathbf{j}_{s+1} = \mathbf{j}_{s}+1$.
 First let us answer the question:

  Are there consecutive jump points of the idempotent endomorphism?

 Yes, for instance $\ep = \wr\,1,1,1,3,5,5\,\wr \in \widehat{\mathcal{E}}_{\mathcal{C}_6}$ is an idempotent and jump points $3$ and $4$ are consecutive. Note that $3$ is also a fixed point of $\ep$. So, we modify the question:

 If the first jump point is not a fixed point, is it possible the next point to be a jump point?

 The answer is negative. Indeed, if $\ep$ is an idempotent and $\ep(\mathbf{j}) = k > j$, then from Proposition 3.1 it follows $\ep(k) = k$ and since $\mathbf{j} + 1 \leq k$ we have $\ep(\mathbf{j}+ 1) \leq k$. So, $\mathbf{j} + 1$ is not a jump point of $\ep$.

\vspace{3mm}

\textbf{Lemma 5.1 } \textsl{ Let $\ep$ be an idempotent and $\mathbf{j}_s$ and $\mathbf{j}_{s+1}$ be non consecutive jump points of $\ep$. Then in the interval $[\mathbf{j}_s,\,\mathbf{j}_{\,s+1} - 1]$ one of the following holds:}

\textsl{1. $\ep$ is a constant endomorphism.}

\textsl{2. $\ep$ is an identity.}

\textsl{3. $\ep$ is an identity in interval $[\mathbf{j}_s,k]$ and a constant endomorphism in interval $[k,\,\mathbf{j}_{\,s+1} - 1]$.}

\textsl{4. $\ep$ is a constant endomorphism in interval $[\mathbf{j}_s,k]$ and an identity in interval $[k,\,\mathbf{j}_{\,s+1} - 1]$.}

\textsl{5. $\ep$ is a constant endomorphism in interval $[\mathbf{j}_s,k]$, an identity in interval $[k,\ell]$ and a constant endomorphism in interval $[\ell,\,\mathbf{j}_{\,s+1} - 1]$.}

\vspace{1mm}

\emph{Proof.} Let $\ep(\mathbf{j}_s) = \mathbf{j}_s$. Since $\mathbf{j}_s + 1$ is not a jump point, there are two possibilities: $\ep(\mathbf{j}_s + 1) = \mathbf{j}_s$ or $\ep(\mathbf{j}_s + 1) = \mathbf{j}_s + 1$. In the first case the equality $\ep(\mathbf{j}_s + 2) = \mathbf{j}_s + 1$ is impossible, see Proposition 3.1. So, either $\ep(\mathbf{j}_s + 2) \geq \mathbf{j}_s + 2$ and $\mathbf{j}_s + 2$ is the next jump point, or $\ep(\mathbf{j}_s + 2) = \mathbf{j}_s$. Using these reasonings, it follows that $\ep$ is a constant endomorphism in  interval $[\mathbf{j}_s,\,\mathbf{j}_{\,s+1} - 1]$. When $\ep(\mathbf{j}_s + 1) = \mathbf{j}_s + 1$ for  point $\mathbf{j}_s +2$ there are three possibilities. If $\ep(\mathbf{j}_s + 2) > \mathbf{j}_s + 2$, then $\mathbf{j}_s + 2$ is the next jump point and $\ep$ is an identity in the interval between considered jump points. If $\ep(\mathbf{j}_s + 2) = \mathbf{j}_s + 1$, using the above reasonings, we find that $\ep$ is an identity in interval $[\mathbf{j}_s,\,\mathbf{j}_s+1]$ and a constant endomorphism in the second interval. Let
 $\ep(\mathbf{j}_s + 2) = \mathbf{j}_s + 2$. Now if for any $k \leq \mathbf{j}_{s+1} - 1$ it follows  $\ep(\mathbf{j}_s + k) = \mathbf{j}_s + k$ then $\ep$ is an identity in the whole interval. If for some $k$ we have $\ep(\mathbf{j}_s + k + 1) = \ep(\mathbf{j}_s + k) = \mathbf{j}_s + k$, then $\ep$ is an identity in interval $[\mathbf{j}_s,k]$ and a constant endomorphism in interval $[k,\,\mathbf{j}_{\,s+1} - 1]$.

Let $\ep(\mathbf{j}_s) = k > \mathbf{j}_s$. Since $\ep$ is an idempotent, it follows $\ep(k) =k$, that is $\ep$ is a constant endomorphism  in interval $[\mathbf{j}_s,k]$. If $\ep(k+1) = k$ then $\ep$ is a constant endomorphism  in the whole interval. Let $\ep(k+1) = k+1$ and for any $\ell$, where
 $\ell = k + 1, \ldots, \mathbf{j}_{s+1} - 1$ we have $\ep(\ell) = \ell$. So, $\ep$ is a constant endomorphism in interval $[\mathbf{j}_s,k]$ and an identity in interval $[k,\,\mathbf{j}_{s+1} - 1]$. In the end if there is some $\ell$, $\ell = k + 1, \ldots, \mathbf{j}_{s+1} - 1$,  such that $\ep(\ell + 1) = \ep(\ell) = \ell$, then $\ep$ is a constant endomorphism in interval $[\mathbf{j}_s,k]$, an identity in interval $[k,\ell]$ and a constant endomorphism in interval $[\ell,\,\mathbf{j}_{s+1} - 1]$.

\vspace{3mm}

It is straightforward to show that if one of the conditions 1. -- 5. of Lemma 5.1 holds for endomorphism $\ep$, then $\ep$ is an idempotent.

\vspace{3mm}

From Lemma 5.1 it follows that the graph of the arbitrary idempotent endomorphism has the form displayed on Figure 1. Note that each of the three segments of this graph may have length zero.

\vspace{1mm}

\begin{figure}[h]\centering
  \includegraphics[width=50mm]{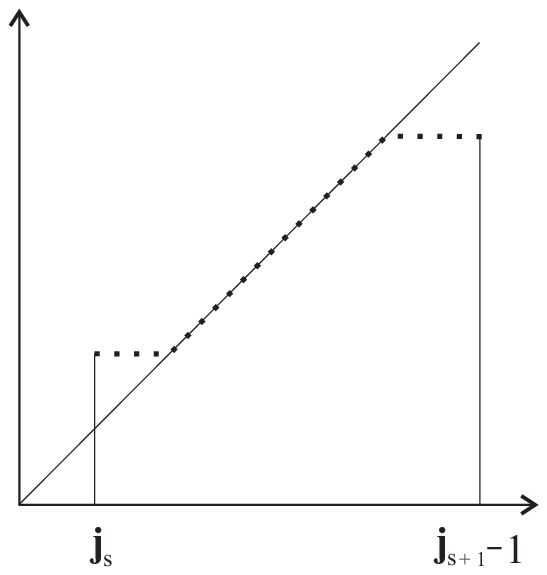}\\
\end{figure}

\vspace{3mm}

\centerline{\small Figure 1 }

\vspace{3mm}

Let $\mathbf{j}_s$ and $\mathbf{j}_{s+1}$ be consecutive jump points of the idempotent $\ep$. Then the graph of $\ep$ is a particular case of those displayed on Figure 1 containing one point $\left(\mathbf{j}_s,\ep(\mathbf{j}_s)\right)$, where $\ep(\mathbf{j}_s) = \mathbf{j}_s$.

\vspace{3mm}

\textbf{Lemma 5.2 } \textsl{Let $\ep$ and $\bar{\ep}$ be idempotent endomorphisms of $\widehat{\mathcal{E}}_{\mathcal{C}_n}$ with the same jump points $\mathbf{j}_1, \ldots, \mathbf{j}_r$. Then $\ep + \bar{\ep}$ is also an idempotent endomorphism with the same jump points.
}

\emph{Proof.} First we shall prove that $\mathbf{j}_1, \ldots, \mathbf{j}_r$ are jump points of the endomorphism $\ep + \bar{\ep}$.
Let $\mathbf{j}$ be any common jump point of $\ep$ and $\bar{\ep}$. If $\mathbf{j}$ is also a fixed point of $\ep$ and $\bar{\ep}$, i.e.  both idempotents satisfied the inequalities $\ep(\mathbf{j} - 1) < \mathbf{j}-1$ and $\ep(\mathbf{j})\geq \mathbf{j}$, $\bar{\ep}(\mathbf{j} - 1) < \mathbf{j}-1$ and $\bar{\ep}(\mathbf{j})\geq \mathbf{j}$, then $(\ep + \bar{\ep})(\mathbf{j} - 1) < \mathbf{j}-1$ and $(\ep + \bar{\ep})(\mathbf{j})\geq \mathbf{j}$. When for $\ep$ it follows $\ep(\mathbf{j} - 1) \leq \mathbf{j}-1$ and $\ep(\mathbf{j}) > \mathbf{j}$ or for $\bar{\ep}$ it follows $\bar{\ep}(\mathbf{j} - 1) \leq \mathbf{j}-1$ and $\bar{\ep}(\mathbf{j})> \mathbf{j}$, then $(\ep + \bar{\ep})(\mathbf{j} - 1) \leq \mathbf{j}-1$ and $(\ep + \bar{\ep})(\mathbf{j}) > \mathbf{j}$.

Now we shall prove that sum $\ep + \bar{\ep}$ is  an idempotent endomorphism. Let us, like in Lemma 5.1, fix  two non consecutive jump points $\mathbf{j}_s$ and $\mathbf{j}_{s+1}$ of idempotents $\ep$ and $\bar{\ep}$ and the interval $[\mathbf{j}_s,\,\mathbf{j}_{s+1} - 1]$.
 When the idempotent $\ep$ satisfies some of conditions 1. -- 5. of Lemma 5.1 we denote this by $\ep_{(1)}, \ldots, \ep_{(5)}$. When $\ep$ satisfies 1. and $\bar{\ep}$ satisfies 5. we denote this by $\ep_{(1)}$ {\small $\&$} $\bar{\ep}_{(5)}$. Let, in this case, $\ep$ be equal to constant $C$ in  interval $[\mathbf{j}_s,\,\mathbf{j}_{s+1} - 1]$ and $\bar{\ep}$ is consecutively equal to constant $k$, identity $\mathbf{i}$ and constant $\ell$. If $C \leq k$ or $C \geq \ell$ then $\ep + \bar{\ep} = C$. If $k < C < \ell$ then $\ep + \bar{\ep}$ satisfy condition 5. This reasoning we denote by $\ep_{(1)}$ {\small $\&$} $\bar{\ep}_{(5)}$ $\Rightarrow (\ep + \bar{\ep})_{(1)}$ or $(\ep + \bar{\ep})_{(5)}$.

  The verification of all fifteen possibilities is easy:

\vspace{2mm}

 $\ep_{(1)} \; \mbox{\small{$\&$}} \; \bar{\ep}_{(1)}\; \Rightarrow\; (\ep + \bar{\ep})_{(1)}$,

 $\ep_{(1)} \; \mbox{\small{$\&$}} \; \bar{\ep}_{(2)}\; \Rightarrow\; (\ep + \bar{\ep})_{(1)} \; \mbox{or}\; (\ep + \bar{\ep})_{(2)} \; \mbox{or} \; (\ep + \bar{\ep})_{(4)}$,

 $\ep_{(1)} \; \mbox{\small{$\&$}} \; \bar{\ep}_{(3)}\; \Rightarrow\; (\ep + \bar{\ep})_{(1)} \; \mbox{or}\; (\ep + \bar{\ep})_{(3)} \; \mbox{or} \; (\ep + \bar{\ep})_{(5)}$,

 $\ep_{(1)} \; \mbox{\small{$\&$}} \; \bar{\ep}_{(4)}\; \Rightarrow\; (\ep + \bar{\ep})_{(1)} \; \mbox{or}\; (\ep + \bar{\ep})_{(4)} \; \mbox{or} \; (\ep + \bar{\ep})_{(5)}$,

 $\ep_{(1)}$ {\small $\&$} $\bar{\ep}_{(5)}$ $\,\Rightarrow\; (\ep + \bar{\ep})_{(1)}$ or $(\ep + \bar{\ep})_{(5)}$, \hspace{2.5mm} $\ep_{(2)} \; \mbox{\small{$\&$}} \; \bar{\ep}_{(2)}\; \Rightarrow\; (\ep + \bar{\ep})_{(2)}$,

 $\ep_{(2)} \; \mbox{\small{$\&$}} \; \bar{\ep}_{(3)}\; \Rightarrow\; (\ep + \bar{\ep})_{(2)}$, \hspace{3mm}  \noindent$\ep_{(2)} \; \mbox{\small{$\&$}} \; \bar{\ep}_{(4)}\; \Rightarrow\; (\ep + \bar{\ep})_{(4)}$, \hspace{3mm}  \noindent$\ep_{(2)} \; \mbox{\small{$\&$}} \; \bar{\ep}_{(5)}\; \Rightarrow\; (\ep + \bar{\ep})_{(4)}$,

 $\ep_{(3)} \; \mbox{\small{$\&$}} \; \bar{\ep}_{(3)}\; \Rightarrow\; (\ep + \bar{\ep})_{(3)}$, \hspace{3mm}  \noindent$\ep_{(3)} \; \mbox{\small{$\&$}} \; \bar{\ep}_{(4)}\; \Rightarrow\; (\ep + \bar{\ep})_{(4)}$, \hspace{3mm}  \noindent$\ep_{(3)} \; \mbox{\small{$\&$}} \; \bar{\ep}_{(5)}\; \Rightarrow\; (\ep + \bar{\ep})_{(5)}$,

 $\ep_{(4)} \; \mbox{\small{$\&$}} \; \bar{\ep}_{(4)}\; \Rightarrow\; (\ep + \bar{\ep})_{(4)}$, \hspace{3mm}  \noindent$\ep_{(4)} \; \mbox{\small{$\&$}} \; \bar{\ep}_{(5)}\; \Rightarrow\; (\ep + \bar{\ep})_{(4)}$, \hspace{3mm}  \noindent$\ep_{(5)} \; \mbox{\small{$\&$}} \; \bar{\ep}_{(5)}\; \Rightarrow\; (\ep + \bar{\ep})_{(5)}$.

\vspace{2mm}

Thus we prove that the sum of two idempotent endomorphisms in interval $[\mathbf{j}_s,\,\mathbf{j}_{s+1} - 1]$ is also idempotent in this interval. So, it follows that $\ep + \bar{\ep}$ is  an idempotent endomorphism in  interval $[\mathbf{j}_{\,1},\,\mathbf{j}_{\,r}]$.
For $0 \leq x < \mathbf{j}_{\,1}$ each of idempotents $\ep$ and $\bar{\ep}$ satisfies one of conditions 1. -- 5. of Lemma 5.1. So, by above reasonings it follows that  sum $\ep + \bar{\ep}$ is an idempotent endomorphism in $[0,\,\mathbf{j}_{\,1} - 1]$.  In a similar way we show that $\ep + \bar{\ep}$ is an idempotent endomorphism in  $[\mathbf{j}_{\,r},\,n-1]$.
 Hence, $\ep + \bar{\ep}$ is an idempotent endomorphism in the whole $\mathcal{C}_n$.

\vspace{3mm}

\textbf{Lemma 5.3 } \textsl{Let $\ep$ and $\bar{\ep}$ be idempotent endomorphisms of $\widehat{\mathcal{E}}_{\mathcal{C}_n}$ with the same jump points $\mathbf{j}_1, \ldots, \mathbf{j}_r$. Then $\ep\cdot \bar{\ep}$ is also an idempotent endomorphism with the same jump points.
}

\emph{Proof.} First we shall prove that $\mathbf{j}_1, \ldots, \mathbf{j}_r$ are the jump points of  endomorphism $\ep\cdot\bar{\ep}$.
Let $\mathbf{j}$ be any common jump point of $\ep$ and $\bar{\ep}$. Let $\mathbf{j}$ be also a fixed point of $\ep$ or of $\bar{\ep}$, that is $\ep(\mathbf{j} - 1) < \mathbf{j}-1$ and $\ep(\mathbf{j})\geq \mathbf{j}$ or $\bar{\ep}(\mathbf{j} - 1) < \mathbf{j}-1$ and $\bar{\ep}(\mathbf{j})\geq \mathbf{j}$. Then it follows in all cases $(\ep\cdot \bar{\ep})(\mathbf{j} - 1) = \bar{\ep}(\ep(\mathbf{j} - 1)) <  \mathbf{j} - 1$ and $(\ep \cdot\bar{\ep})(\mathbf{j}) = \bar{\ep}(\ep(\mathbf{j})) \geq \mathbf{j}$. Let $\mathbf{j}$ be not a fixed point of both idempotents, i.e
$\ep(\mathbf{j} - 1) \leq \mathbf{j}-1$ and $\ep(\mathbf{j}) > \mathbf{j}$ and also $\bar{\ep}(\mathbf{j} - 1) \leq \mathbf{j}-1$ and $\bar{\ep}(\mathbf{j})> \mathbf{j}$. At last we obtain that $(\ep\cdot \bar{\ep})(\mathbf{j} - 1) = \bar{\ep}(\ep(\mathbf{j}-1)) \leq \bar{\ep}(\mathbf{j}-1)) \leq \ep(\mathbf{j}-1)$ and $(\ep\cdot \bar{\ep})(\mathbf{j}) = \bar{\ep}(\ep(\mathbf{j} - 1)) > \bar{\ep}(\mathbf{j}-1) > \ep(\mathbf{j}-1)$.

Now we shall prove that  $\ep\cdot \bar{\ep}$ is  an idempotent endomorphism. Like in the proof of the previous lemma, we  fix  two non consecutive jump points $\mathbf{j}_s$ and $\mathbf{j}_{s+1}$ of idempotents $\ep$ and $\bar{\ep}$ and  interval $[\mathbf{j}_s,\,\mathbf{j}_{s+1} - 1]$.

Using the notations from the proof of Lemma 5.2 we easily verify all twenty five possibilities:
$$\begin{array}{llll} \ep_{(1)}\cdot \bar{\ep}_{(1)} = \ep_{(1)}, & \ep_{(1)}\cdot \bar{\ep}_{(2)} = \ep_{(1)}, & \ep_{(1)} \; \mbox{\small{$\&$}} \; \bar{\ep}_{(3)}\; \Rightarrow\; (\ep\cdot\bar{\ep})_{(1)}, & \ep_{(1)} \; \mbox{\small{$\&$}} \; \bar{\ep}_{(4)}\; \Rightarrow\; (\ep\cdot\bar{\ep})_{(1)},\\
\ep_{(1)} \; \mbox{\small{$\&$}} \; \bar{\ep}_{(5)}\; \Rightarrow\; (\ep\cdot\bar{\ep})_{(1)}, & \ep_{(2)}\cdot \bar{\ep}_{(1)} = \bar{\ep}_{(1)}, & \ep_{(2)}\cdot \bar{\ep}_{(2)} = \bar{\ep}_{(2)}, & \ep_{(2)}\cdot \bar{\ep}_{(3)} = \bar{\ep}_{(3)},\\
\ep_{(2)}\cdot \bar{\ep}_{(4)} = \bar{\ep}_{(4)}, & \ep_{(2)}\cdot \bar{\ep}_{(5)} = \bar{\ep}_{(5)}, & \ep_{(3)}\cdot \bar{\ep}_{(1)} = \bar{\ep}_{(1)}, & \ep_{(3)}\cdot \bar{\ep}_{(2)} = \ep_{(3)},\\
\end{array} $$

$$\begin{array}{ll} \ep_{(3)} \; \mbox{\small{$\&$}} \; \bar{\ep}_{(3)}\; \Rightarrow\; (\ep\cdot\bar{\ep})_{(3)}, & \ep_{(3)} \; \mbox{\small{$\&$}} \; \bar{\ep}_{(4)}\; \Rightarrow\; (\ep\cdot\bar{\ep})_{(1)} \; \mbox{or}\; (\ep\cdot\bar{\ep})_{(5)},\\  \ep_{(3)} \; \mbox{\small{$\&$}} \; \bar{\ep}_{(5)}\; \Rightarrow\; (\ep\cdot\bar{\ep})_{(1)} \; \mbox{or}\; (\ep\cdot\bar{\ep})_{(5)}, & \ep_{(4)}\cdot \bar{\ep}_{(1)} = \bar{\ep}_{(1)},\\
\ep_{(4)}\cdot \bar{\ep}_{(2)} = \ep_{(4)}, & \ep_{(4)} \; \mbox{\small{$\&$}} \; \bar{\ep}_{(3)}\; \Rightarrow\; (\ep\cdot\bar{\ep})_{(1)} \; \mbox{or}\; (\ep\cdot\bar{\ep})_{(5)},\\
 \ep_{(4)} \; \mbox{\small{$\&$}} \; \bar{\ep}_{(4)}\; \Rightarrow\; (\ep\cdot\bar{\ep})_{(4)}, & \ep_{(4)} \; \mbox{\small{$\&$}} \; \bar{\ep}_{(5)}\; \Rightarrow\; (\ep\cdot\bar{\ep})_{(1)} \; \mbox{or}\; (\ep\cdot\bar{\ep})_{(5)},\\
\end{array}$$

$$\begin{array}{lll} \ep_{(5)}\cdot \bar{\ep}_{(1)} = \bar{\ep}_{(1)}, &\;\;\, \ep_{(5)}\cdot \bar{\ep}_{(2)} = \ep_{(5)}, & \ep_{(5)} \; \mbox{\small{$\&$}} \; \bar{\ep}_{(3)}\; \Rightarrow\; (\ep\cdot\bar{\ep})_{(1)} \; \mbox{or}\; (\ep\cdot\bar{\ep})_{(5)},\\
\end{array}$$
$$\begin{array}{ll} \ep_{(5)} \; \mbox{\small{$\&$}} \; \bar{\ep}_{(4)}\; \Rightarrow\; (\ep\cdot\bar{\ep})_{(1)} \; \mbox{or}\; (\ep\cdot\bar{\ep})_{(5)}, & \ep_{(5)} \; \mbox{\small{$\&$}} \; \bar{\ep}_{(5)}\; \Rightarrow\; (\ep\cdot\bar{\ep})_{(1)} \; \mbox{or}\; (\ep\cdot\bar{\ep})_{(5)}.\\
\end{array}$$

Thus we prove that the product of two idempotent endomorphisms in interval $[\mathbf{j}_s,\,\mathbf{j}_{s+1} - 1]$ is also idempotent in this interval. It easily implies  that $\ep\cdot \bar{\ep}$ is  an idempotent endomorphism in the whole interval $[\mathbf{j}_{\,1},\,\mathbf{j}_{\,r}]$.
In  interval $[0,\,\mathbf{j}_{\,1} - 1]$ each of idempotents $\ep$ and $\bar{\ep}$ satisfies one of conditions 1. -- 5. of Lemma 5.1. So, by above reasonings it follows that the product $\ep\cdot \bar{\ep}$ is an idempotent endomorphism in $[0,\,\mathbf{j}_{\,1} - 1]$.  In a similar way we show that $\ep\cdot\bar{\ep}$ is an idempotent endomorphism in  $[\mathbf{j}_{\,r},\,n-1]$.
 Hence, $\ep\cdot \bar{\ep}$ is an idempotent endomorphism in the whole $\mathcal{C}_n$.

\vspace{3mm}

Immediately from Lemma 5.2 and Lemma 5.3 it follows

\vspace{3mm}

\textbf{Theorem 5.4} \textsl{The set of idempotent endomorphisms of $\widehat{\mathcal{E}}_{\mathcal{C}_n}$ with the same jump points is a subsemiring of $\widehat{\mathcal{E}}_{\mathcal{C}_n}$.
}

\vspace{3mm}

Using Lemma 5.2, Lemma 5.3 and Figure 1 it is easy to prove

\vspace{3mm}

\textbf{Proposition 5.5} \textsl{The set of idempotent endomorphisms of $\widehat{\mathcal{E}}_{\mathcal{C}_n}$ without jump points is a subsemiring of $\widehat{\mathcal{E}}_{\mathcal{C}_n}$ of order $\ds \binom{n+1}{2}$}.

\vspace{8mm}

\centerline{\large References}

\vspace{4mm}

[1] J. Almeida, S. Margolis, B. Steinberg and  M. Volkov, Representation theory of finite semigroups, semigroup radicals and formal
language theory, Trans. Amer. Math. Soc. 361 (3) (2009) 1429 -- 1461.

[2]  O. Ganyushkin and V. Mazorchuk, Classical Finite Transformation Semigroups: An Introduction,
Springer-Verlag London Limited, 2009.

[3]  J. Golan, Semirings and Their Applications, Kluwer, Dordrecht, 1999.

[4] Z. Izhakian, J. Rhodes and  B. Steinberg, Representation theory of finite semigroups over semirings, Journal of Algebra 336 (2011) 139 -- 157.

[5]  J. Je$\hat{\mbox{z}}$ek, T. Kepka and  M. Mar\`{o}ti, The endomorphism semiring of a se\-milattice,
Semigroup Forum, 78 (2009), 21 -- 26.

[6]  E. H.  Moore, A definition of abstract groups, Trans. Amer. Math. Soc, 3 (1902), 485 -- 492.

[7] J. Rhodes and B. Steinberg, The q-Theory of Finite Semigroups, Springer Monogr. Math., Springer, New York, 2009.

[8] R. Stanley, Enumerative combinatorics, Vol. 2, Cambridge University Press, 1999.

[9]  I. Trendafilov and  D. Vladeva, The endomorphism semiring of a finite chain, Proc.
Techn. Univ.-Sofia, 61, 1, (2011), 9 -- 18.

[10]  I. Trendafilov and  D. Vladeva, Subsemirings of the endomorphism semiring of a finite chain, Proc.
Techn. Univ.-Sofia, 61, 1, (2011), 19 -- 28.

[11]  I. Trendafilov and  D. Vladeva, On some semigroups of the partial transformation semigroup, Appl. Math. in Eng. and Econ. -- 38th Int. Conf. (2012) AIP Conf. Proc. (to appear).

[12]  J. Zumbr\"{a}gel, Classification of finite congruence-simple semirings with zero,
J. Algebra Appl. 7 (2008) 363 -- 377.

\vspace{7mm}

 Authors:  
 
 Ivan Trendafilov,  Department "Algebra and geometry"$\hphantom{}$, FAMI, TU--Sofia, \emph{e-mail:}
ivan$\_$d$\_$trendafilov@abv.bg

 Dimitrinka Vladeva,  Department "Mathema\-tics and physics"$\hphantom{}$, LTU, Sofia, \emph{e-mail:}
d$\_$vladeva@abv.bg

\end{document}